\documentclass[12pt,reqno]{amsart}
\usepackage{amsmath,amssymb,amsfonts,amscd,latexsym,amsthm,mathrsfs}
\usepackage[unicode]{hyperref}
\newtheorem{theorem}{Theorem}

\let\wt\widetilde
\begin{document}

\title{A \protect\lowercase{strange identity of an} MF (M\protect\lowercase{ahler function})}

\dedicatory{In memoriam Peter Bundschuh (1938--2024), ein Lehrer und Freund}

\date{20 March 2024}

\author{W\protect\lowercase{adim} Z\protect\lowercase{udilin}}
\address{Department of Mathematics, IMAPP, Radboud University, PO Box 9010, 6500~GL Nijmegen, Netherlands}
\email{w.zudilin@math.ru.nl}


\begin{abstract}
We relate two different solutions of a Mahler equation; one solution is only defined at certain roots of unity, while the other is an analytic function inside the unit disk.
\end{abstract}

\maketitle

As a student and young researcher working in transcendental number theory I was educated in the Siegel--Shidlovsky and Gelfand--Schneider methods, and was aware of the Mahler method but never really paid a serious attention to the latter.
The things changed in 2003 during my stay at the University of Cologne as a Humboldt fellow.
Peter Bundschuh was my host and mentor there, and chatting with him significantly extended my `transcendental' horizons.
There was a big deal of Mahler's legacy in our discussions, in particular, about the Mahler method \cite{Ma29,Ma69,Ni97} and contributions of his own and of his students; Bundschuh's latest works on the topic are \cite{BV15a,BV15b,BV16a,BV16b,BV18a,BV18b,BV18c,BV20} jointly with Keijo V\"a\"an\"anen.
Applications of the method were quite scarce, mainly because there were not so many `natural' examples of what we now call \emph{Mahler functional equations} and \emph{Mahler functions}; this has changed a lot in the last decades \cite{Ad19,Ma19,Ni21}.
One principal question about the Mahler functions is their linear and algebraic independence, as then it immediately implies related results for their values.
It happened much later that I could contribute to this myself, as part of Novocastrian `CARMA' team \cite{BCZ16};
with Richard Brent and Michael Coons we have developed a new machinery for investigating independence by looking at the asymptotics of Mahler functions at roots of unity.
This methodology is classical in the context of modular forms and their numerous generalisations.
My goal here is to pay tribute to the late Peter Bundschuh, through giving another illustration of the similarity between the two quite different functional classes.

\section{A Mahler function}
\label{sec1}

The function
\[
U(q)=\sum_{n=0}^\infty\prod_{j=0}^{n-1}(q^{2^j}-1)
=\sum_{n=0}^\infty(-1)^n\prod_{j=0}^{n-1}(1-q^{2^j})
\]
is a solution of the Mahler functional equation
\begin{equation}
U(q)=1+(q-1)U(q^2).
\label{eq2}
\end{equation}
The function is well defined on the set $\Xi=\{\exp(2\pi\mathrm{i}j/2^m):j\in\mathbb Z,\ m\in\mathbb Z_{\ge0}\}$ of roots of unity but nowhere else. It can be also expanded into a formal power series in $\mathbb Z[[q-1]]$ (or in $\mathbb Z[[q-\xi]]$ for any $\xi\in\Xi$):
\begin{align}
U(q)
&= 1 + (q-1) + 2(q-1)^2 + 9(q-1)^3 + 80(q-1)^4 + 1390(q-1)^5
\nonumber\\ &\quad
+ 47094(q-1)^6 + 3127145(q-1)^7 + 409428448(q-1)^8 + O\big((q-1)^9\big).
\label{eq2a}
\end{align}
At the same time, the Mahler equation \eqref{eq2} possesses a unique power-series solution
\[
U_0(q)= \frac12\big(1 + q - q^2 + q^3 + q^4 - q^5 - q^6 + q^7 - q^8 + q^9 + q^{10} - q^{11} + O(q^{12})\big)
\]
at the origin, so that $2U_0(q)$ represents \emph{all} powers of $q$ and is a $\{\pm1\}$-series.
With a little search one finds out that
$$
U_0(q)+\frac1{2(1-q)}
= 1 + q + q^3 + q^4 + q^7 + q^9 + q^{10} + q^{12} + q^{15} + q^{16} + q^{19} + O(q^{21})
$$
is the generating function of `natural numbers having an even number of nonleading zeros in their binary expansion'
(\href{http://oeis.org/A059010}{OEIS A059010}). It means that the parity of the number of those zeroes contributes to the sign of original $U_0(q)$.
This type of combinatorial interpretation is customary for solutions of Mahler functional equations.

Our principal result is the following `strange identity', in the spirit of Zagier's identity in~\cite{Za01}.

\begin{theorem}
\label{th1}
For the derivatives of all orders, the function $U(q)$ agrees at every root of unity in $\Xi=\{\exp(2\pi\mathrm{i}j/2^m):j\in\mathbb Z,\ m\in\mathbb Z_{\ge0}\}$ with the radial limit of the function $U_0(q)$.
\end{theorem}

\section{Proof}
\label{sec2}

For convenience, introduce the notation
\[
\{q\}_n=\prod_{j=0}^{n-1}(1-q^{2^j})
\]
for the Mahler factorial. This is defined for $n=0,1,2,\dots$ but also for $n=\infty$ when $|q|<1$ or as a power series in~$q$.

Let us consider the partial sums
\[
U(q;N)=\sum_{n=0}^{N-1}(-1)^n\{q\}_n
\]
of $U(q)$. A simple inductive argument shows that
\begin{align*}
U(q;2k+1)
&=1-q^2(1-q)-q^8(1-q)(1-q^2)(1-q^4)-\dotsb-q^{2^{2k-1}}\{q\}_{2k-1}
\\
&=1-\sum_{i=0}^{k-1}q^{2^{2i+1}}\{q\}_{2i+1}
\\ \intertext{and}
U(q;2k)
&=q+q^4(1-q)(1-q^2)+q^{16}\{q\}_4+\dots+q^{2^{2k-2}}\{q\}_{2k-2}
\\
&=\sum_{i=0}^{k-1}q^{2^{2i}}\{q\}_{2i}
\end{align*}
for $k=0,1,2,\dots$\,. Since the degree of $q^{2^i}\{q\}_i$ is equal to $2^{i+1}-1$, the terms in the sum do not mix with each other and the $q$-series
\[
U_-(q)=\lim_{k\to\infty}U(q;2k+1)=1-\sum_{i=0}^\infty q^{2^{2i+1}}\{q\}_{2i+1}
\]
and
\[
U_+(q)=\lim_{k\to\infty}U(q;2k)=\sum_{i=0}^\infty q^{2^{2i}}\{q\}_{2i}
\]
are well defined, and their coefficients belong to the set $\{0,\pm1\}$.
In particular, the series converge in the disk $|q|<1$. They are also well defined on the set $\Xi$ of roots of unity, and their values coincide with those of $U(q)$ on the set.
Further induction shows that the polynomial
\[
\wt U(q;N)=U(q;N)+\frac{(-1)^N}2\{q\}_N
\]
satisfies the functional equation
\[
1+(q-1)\wt U(q^2;N)=\wt U(q;N+1) \quad\text{for}\; N=0,1,2,\dots\,.
\]
By taking the limit as $N\to\infty$ separately along the even and odd $N$ we obtain the following identities for the formal series (hence when $|q|<1$):
\[
U_0(q)=U_+(q)+\frac12\{q\}_\infty=U_-(q)-\frac12\{q\}_\infty.
\]
The theorem now follows from noticing that $\{q\}_\infty|=0$ and recalling that $U_+(q)=U_-(q)=U(q)$ at every $q=\xi\in\Xi$, also for the derivatives of all orders.

\section{Concluding remarks}
\label{sec3}

Unlike the situation in \cite{Za01}, it is hard to expect a combinatorial meaningfulness of the (provably!) positive coefficients in the expansion \eqref{eq2a} of $U(q)$ in powers of $q-1$.
Is their asymptotics interesting?

Our function $U(q)$ is a member of the family $W(q)=W(q;a)=\sum_{n=0}^\infty a^n\{q\}_n$ defined on $\Xi$ and satisfying the Mahler equation $W(q)=1+a(1-q)W(q^2)$.
Assuming that $a\ne1$, the equation also possesses the solution
$$
W_0(q)=W_0(q;a)=\frac1{1-a}\sum_{m=0}^\infty(-1)^{s_2(m)}a^{r_2(m)}q^m,
$$
where $r_2(m)=\big\lfloor\frac{\log m}{\log2}\big\rfloor+1$ is the number of binary digits of $m$ and $s_2(m)$ is their sum (conventionally, $r_2(0)=s_2(0)=0$).
One may notice that no analytic solution to this Mahler equation exists at $\infty$.
There are however many cases when power-series solutions to Mahler equations exist both at $0$ and at $\infty$, and they are strangely (dis)connected at roots of unity; one can find a related discussion, together with new beautiful links to continued fractions, in the thesis \cite{Ni21} of Joris Nieuwveld (see also the commentary in~\cite{Zu21}).

Many other remarkable features of Mahler equations and functions wait to be discovered.


\end{document}